\documentstyle[twoside]{article}
\newtheorem{theorem}{Theorem}[section]

\newtheorem{definition}{Definition}[section]
\newtheorem{proposition}{Proposition}[section]
\newtheorem{corollary}{Corollary}[section]
\newtheorem{lemma}{Lemma}[section]

\newcommand{\proof}{\noindent{\sc Proof. }}
\newcommand{\bbfR}{{I\!\!R}}
\newcommand{\bbfN}{{I\!\!N}}
\renewcommand{\P}{{I\!\!P}}

\renewcommand{\div}{\nabla\cdot}
\renewcommand{\S}{{\cal S}}
\newcommand{\C}{{\cal C}}
\newcommand{\Y}{{\cal Y}}

\newcommand{\rf}[1]{(\ref{#1})}

\newcommand{\X}{{\cal X}}
\newcommand{\Ep}{{E_p}}
\newcommand{\BEalf}{{BE^\alpha}}
\newcommand{\BEpa}{{BE^\alpha_p}}
\newcommand{\BEpn}{{BE_p^{1-3/p}}}

\def\cbdu{\hfill{$\Box$}}

\def\esssup{\mbox{\rm ess} \sup_{\!\!\!\!\!\!\!\!\! \xi\in\bbfR^3}}
\newcommand{\PM}{{{\cal P\! M}}}

\newcounter{remark}
\newenvironment{remark}%
{\medskip\stepcounter{remark}\noindent {\it Remark
\arabic{section}.\arabic{remark}} }{\rm \cbdu\medskip}

\title{\bf About the regularized Navier--Stokes equations }

\author{{\sc Marco Cannone \& Grzegorz Karch}\\
\\
{\small Universit\'e de Marne-la-Vall\'ee, Laboratoire d'Analyse}\\
{\small et de Math\'ematiques Appliqu\'ees}\\
{\small Cit\'e Descartes--5, bd Descartes, Champs-sur-Marne,}\\
{\small 77454 Marne-la-Vall\'ee Cedex 2, France}\\
{\small {\tt cannone@math.univ-mlv.fr}}\\
     \\
{\small Instytut Matematyczny, Uniwersytet Wroc\l awski}\\
{\small pl. Grunwaldzki 2/4, 50-384 Wroc\l aw, Poland;}\\
{\small Institute of Mathematics}\\
{\small Polish Academy of Sciences, Warsaw
(2002-2003)}\\
{\small {\tt karch@math.uni.wroc.pl}}\\
}

\date{}

\begin{document}
\maketitle

\begin{abstract}
The first goal of this paper is to study the large time behavior of 
solutions to the Cauchy problem for the 3-dimensional incompressible 
Navier-Stokes system. The Marcinkiewicz space $L^{3,\infty}$ is used 
to prove some asymptotic stability results for solutions with 
infinite energy.
Next, this approach is applied to the analysis of two classical 
``regularized''  Navier-Stokes systems. The first one was
introduced by J. Leray and consists in ``mollifying'' the 
nonlinearity. The second one was proposed by J.L. Lions, who added
the artificial hyper-viscosity $(-\Delta)^{\ell/2}$, $\ell>2$, to the 
model. It is shown in the present paper that, in the
whole space, solutions to those modified models converge as 
$t\to\infty$ toward solutions of the original Navier-Stokes
system.
\footnote[0]{ 2000 {\it Mathematics Subject Classification}:
35Q30, 76D05, 35B40.
}
\footnote[0]{ {\it Key words and phrases}:
incompressible Navier--Stokes system;
self-similar solutions; the Mar\-cin\-kie\-wicz spaces;  asymptotic stability.}

\end{abstract}

\newpage

\section{Introduction}

Since the seminal paper by Leray \cite{Ler}, several methods have been
developed
to prove existence of global-in-time   weak solutions of the Cauchy
problem for the
three-di\-men\-sion\-al
Navier--Stokes system
\begin{eqnarray}
&&u_t-\Delta u +(u\cdot \nabla) u +\nabla p =F,\quad x\in \bbfR^3,
t>0\label{eq-0}\\
&&\div u =0, \label{div-0}\\
&&u(0)=u_0.\label{ini-0}
\end{eqnarray}
The usual tool is to consider a regularized problem (containing a
parameter $\kappa>0$) for which one can prove the existence of a
unique  smooth solution. Next, due to the energy
inequality, one can pass to the limit as $\kappa\to 0$ and to show
that the limit function is a~weak solution to problem
\rf{eq}--\rf{ini}.
A detailed description of the possible ways used in the literature 
for modifying the system \rf{eq-0}--\rf{ini-0}
is contained in \cite{FrPa}.

This idea was used already by  Leray \cite{Ler},
who mollified equation \rf{eq-0} replacing the nonlinearity 
$\nabla\cdot (u\otimes u)$
by the smoother term $\nabla\cdot ((u*\omega_\kappa)\otimes u)$ with a
smooth function $\omega$ such that $\int \omega \;dx =1$
and $\omega_\kappa (x)=\kappa^{-3} \omega(x/\kappa)$.
On the other hand, J.-L. Lions proposed to replace the Laplacian
$-\Delta$  by the sum $-\Delta +\kappa (-\Delta)^{\ell/2}$,  $\ell>
2$ (in a
way that is
reminiscent of a Taylor expansion), and
for such a modified problem considered in a bounded domain,
J.-L. Lions was able to prove (cf. \cite[Chap. 1, Remarque
6.11]{Lio3}) the existence of a unique regular solution provided
$\ell \geq 5/2$ ($\ell\geq (n+2)/2$ for the $n$-dimensional
problem).
An analogous result for the whole space $\bbfR^3$ is contained e.g. 
in \cite{KP}.
Hence, one can say that the mollified nonlinearity as well as the 
hyperdissipative
term  in the equation smooth out solutions.

The goal of this paper is to show that,
in the whole space $\bbfR^3$, such  corrections in the model
disappear
asymptotically as $t\to\infty$, at least, when small solutions are
considered.
More precisely, we fix $\kappa>0$ in both models, and we show that their
solutions converge in a suitable sense as $t\to\infty$ toward solutions of
the Navier-Stokes system \rf{eq-0}--\rf{ini-0} corresponding to the same
initial conditions and external forces.

\medskip
{\bf Notations.}
The notations to be used are mostly standard. For
$1\leq p\leq \infty$, the $L^p$-norm   of a Lebesgue
measurable real-valued function defined on $\bbfR^3$
is denoted by $\|v\|_p$. On the other  hand, the norm of the weak $L^p$-
space (the Marcinkiewicz space) $L^{p,\infty}=L^{p,\infty}(\bbfR^3)$ 
is denoted by $\|\cdot\|_{p,\infty}$;
cf.~Section 3 for suitable definitions.
We will always  denote by
$\|\cdot\|_{\cal X}$ the norm of any other Banach space $\cal X$
used in this paper. Here, we study properties of  vector-valued solutions
$u=(u_1,u_2,u_3)$ to the Navier-Stokes system \rf{eq-0}--\rf{ini-0}, hence
the notation $u\in \X$ should be understood as $u_i\in \X$ for every
$i=1,2,3$; moreover, by the very definition, $\|u\|_\X =\max\{\|u_1\|_\X,
\|u_2\|_\X, \|u_3\|_\X\}$.


\section{Results and comments}
\setcounter{equation}{0}
\setcounter{remark}{0}

Let us recall  the projection
$\P$
of $(L^2)^3$ onto the
subspace $\P [(L^2)^3]$
of  solenoidal vector fields
(i.e. those characterized by the divergence
condition \rf{div-0}). It is known that $\P$ is a pseudodifferential
operator of order 0. In fact, it can be written
     as a~combination of the Riesz transforms $R_j$
with symbols $\xi_j/|\xi|$,
$$
\P(v_1, v_2, v_3)=(v_1-R_1\sigma, v_2-R_2\sigma, v_3-R_3\sigma),
$$
where $\sigma = R_1v_1+R_2v_2+R_3v_3$. This explicit formula
allows us to consider $\P$ as the bounded operator on $L^p=L^p(\bbfR^3)$, $1<p<\infty$, as
well
as on the Marcinkiewicz weak $L^p$-spaces recalled in the next section.

Using this projection, one can remove the pressure from the model
\rf{eq-0}-\rf{ini-0} and  obtain an equivalent Cauchy problem
\begin{eqnarray}
&&u_t-\Delta u +\P\nabla\cdot (u\otimes u)  =\P F,\quad x\in \bbfR^3,
t>0\label{eq}\\
&&u(0)=u_0.\label{ini}
\end{eqnarray}

Our first goal is to study  solutions
to  problem \rf{eq}--\rf{ini} rewritten as
  the integral equation
\begin{eqnarray}
u(t)&=& S(t) u_0 - \int_0^t  S(t-\tau) \P \nabla \cdot (u\otimes
u)(\tau)\;d\tau \label{duhamel}\\
& & +\; \int_0^t S(t-\tau) \P F(\tau)\;d\tau. \nonumber
\end{eqnarray}
Here, the heat semigroup on $\bbfR^3$, denoted by $S(t)$, is realized as the
convolution with the Gaussian kernel $p(x,t)=(4\pi t)^{-3/2}
\exp(-|x|^2/(4t))$.
Note that \rf{duhamel} has the form $u=y+B(u,u)$, where the bilinear form
is defined
as
\begin{equation}
B(u,v)(t)=-\int_0^t S(t-\tau) \P \nabla \cdot (u\otimes
v)(\tau)\;d\tau,\label{Bf}
\end{equation}
and $y=S(t) u_0 +\int_0^t S(t-\tau) \P F(\tau)\;d\tau$.
Hence,  using  the classical Picard approach,
which is based on  Lemma \ref{lem:xyB} below, one can easily construct
solutions
in the space
\begin{equation}
\X_3=C_w([0, \infty), L^{3,\infty})\label{X3}
\end{equation}
provided initial data and external forces are small in a suitable sense.

Here, it should be emphasized  that
the Marcinkiewicz space $L^{p,\infty}$ is not separable and the heat
semigroup is not
strongly continuous on the space. Hence, in our considerations below,
  we introduce
the space
$C_w([0,\infty), \;L^{p,\infty})$ consisting of functions $u$ with the
following two properties
\begin{itemize}
\item $u$ is bounded and continuous from $(0,\infty)$ to $L^{p,\infty}$ in the norm
topology of $L^{p,\infty}$;
\item $u(\cdot,t)\to u(\cdot,0)$ weakly as $t\searrow 0$.
\end{itemize}
We refer the reader e.g. to \cite[Lemmata 2 \& 3]{Ba1} for the detailed  proof that $S(\cdot)u_0\in
C_w([0,\infty), \;L^{p,\infty})$ for every $u_0\in L^{p,\infty}$.

To deal with external forces and for simplicity of notation we 
introduce another space
\begin{eqnarray}
\Y_3 =\bigg\{F\in C([0,\infty), \S'(\bbfR^3))&:& \mbox{the integral} 
\;\; \int_0^t S(t-\tau)\P F(\tau)\;d\tau \label{Y3}\\
&& \mbox{makes sense and belongs to}\;\; \X_3\bigg\} \nonumber
\end{eqnarray}
supplemented with the norm
$$
\|F\|_{\Y_3}\equiv \sup_{t>0} \left\|\int_0^t S(t-\tau)\P 
F(\tau)\;d\tau\right\|_{3,\infty}=
\left\|\int_0^\cdot S(\cdot-\tau)\P F(\tau)\;d\tau\right\|_{\X_3}.
$$

Below, at the end of Section 3, we discuss sufficient conditions 
which guarantee that $F\in \Y_3$ thus showing that
$\Y_3\not = \emptyset $.

Now, the result on the existence of small solutions in the space 
$\X_3$ reads as follows.

\begin{theorem}\label{th:eq}
Assume that $u_0\in L^{3,\infty}$ and $F\in \Y_3$
satisfy $\|u_0\|_{3,\infty}+\|F\|_{\Y_3}<\varepsilon$
for some $0<\varepsilon<1/(4\eta)$ where the constant $\eta$ is defined in
Proposition \ref{prop:L-S}, below.
There exists a global-in-time solution of
\rf{eq}--\rf{ini}  in the space
$
\X_3.
$
This is the unique solution satisfying the condition
$ \|u\|_{\X_3}\leq 2 \varepsilon.$
Moreover, this solution depends continuously on initial data
and external forces in the following sense: for every two solutions 
$u$ and $\tilde u$ of \rf{eq}-\rf{ini} corresponding to initial data 
$u_0$, $\tilde u_0$ and forces $F$, $\tilde F$, respectively, 
satisfying the assumptions of the first part of this theorem, we have
$$
\|u-\tilde u\|_{\X_3}\leq C(\|u_0-\tilde u_0\|_{3,\infty}+\|F-\tilde 
F\|_{\Y_3} )
$$
with a constant $C=C(\varepsilon, \eta)=1/(1-4\eta\varepsilon)$.
\cbdu
\label{th:glob}
\end{theorem}

The regularity of small solutions depends essentially on the external
forces. Indeed, in Proposition \ref{eq:reg} below,
  we show that solutions constructed in Theorem \ref{th:eq}  belong, 
moreover, to the space
\begin{equation}
\X_p =\{ u\in C((0,\infty), L^{p,\infty}(\bbfR^3))\;:\;
\sup_{t>0} t^{(1-3/p)/2} \|u(t)\|_{p,\infty}<\infty\}
\label{Xp}
\end{equation}
provided  $\int_0^\cdot S(\cdot-\tau)\P F(\tau)\;d\tau\in \X_3\cap \X_p$.
Here, we obtain immediately a better regularity result due to the
  imbedding
\begin{equation}
  L^{3,\infty}\cap L^{p,\infty}\subset L^q,\quad \mbox{for every}\quad  3<q<p.
\label{imbedding}
\end{equation}
Since, by Proposition \ref{eq:reg}, $u\in \X_3\cap \X_p$, we easily 
deduce from \rf{imbedding}
the decay rates of solutions in the Lebesgue space $L^q$:
$$
\sup_{t>0} t^{(1-3/q)/2} \|u(t)\|_q<\infty, \quad \mbox{for 
every}\quad  q\in (3,\infty).
$$

The following theorem is the new contribution to the theory concerning
large time behavior of solutions discussed above.

\begin{theorem}\label{asymp}
Let the assumptions of Theorem \ref{th:glob} hold true.
Assume that $u$ and $\tilde u$ are two solutions of \rf{eq}--\rf{ini}
constructed in Theorem \ref{th:glob}
corresponding to the initial conditions $u_0,\tilde u_0\in
L^{3,\infty}$ and external forces $F,\tilde F\in \Y_3$, respectively.
Suppose that
\begin{equation}
\lim_{t\to\infty}\left\|S(t)(u_0-\tilde u_0)
+\int_0^t S(t-\tau)\P (F(\tau)-\tilde F(\tau))\;d\tau\right\|_{3,\infty}=0.
\label{e:u0v0}
\end{equation}
Then
\begin{equation}
\lim_{t\to\infty}\|u(\cdot,t)-\tilde u(\cdot,
t)\|_{3,\infty}=0\label{e:utvt}
\end{equation}
holds.
\end{theorem}

As will be proved in Corollary~\ref{or:as} below, conditions 
\rf{e:u0v0} and \rf{e:utvt} are, in fact, equivalent.

Section 4 contains more results being direct corollaries of Theorem
\ref{asymp}. In particular, it is shown that under the assumptions of this
theorem
\begin{equation}
\lim_{t\to\infty}t^{(1-3/p)/2}\|u(\cdot,t)-\tilde u(\cdot,t)\|_p=0\label{Lp-as}
\end{equation}
for every $p\in (3,\infty)$. First, we show relation \rf{Lp-as}
with the Lebesgue norm replaced by
  the Mar\-cin\-kie\-wicz
$L^{p,\infty}$-norm, next,
the limit in \rf{Lp-as}
  results directly from the imbedding \rf{imbedding}.
Some details are contained in Proposition \ref{prop:Lp}, below.

\medskip

Theorem \ref{asymp} is the counterpart of a result contained in 
\cite{CK1} where global-in-time solutions are
constructed in the space
$$
\PM^2\equiv\{v\in{\cal S}'(\bbfR^d): \widehat v\in L^1_{\rm loc}(\bbfR^d),
\|v\|_{\PM^2}\equiv \esssup |\xi|^2|\widehat
v(\xi)|<\infty\}.
$$
In particular, in that setting, it is possible to study one-point 
stationary singular solutions to \rf{eq}--\rf{ini}
(constructed independently by Landau and Tian and Xin) of the 
following form (cf. \cite{TX})
\begin{eqnarray}
u_1(x)= 2{c|x|^2-2x_1|x|+cx_1^2\over |x|(c|x|-x_1)^2},  &&\quad
u_2(x)= 2{x_2(cx_1-|x|) \over |x|(c|x|-x_1)^2}, \label{sing-sol} \\
u_3(x)= 2{x_3(cx_1-|x|)  \over |x|(c|x|-x_1)^2},&&\quad
p(x)= 4{cx_1-|x|  \over |x|(c|x|-x_1)^2} \nonumber
\end{eqnarray}
where $|x|=\sqrt{x_1^2+x_2^2+x_3^2}$ and $c$ is an arbitrary constant
such that $|c|>1$.
By straightforward calculations, one can check that,
indeed, the functions
$u_1(x)$, $u_2(x)$, $u_3(x)$, and $p(x)$ given by \rf{sing-sol} satisfy
\rf{eq}--\rf{ini} with
$F\equiv 0$ in the {\sl pointwise sense} for every $x\in
\bbfR^3\setminus
\{(0,0,0)\}$. On the other hand, if one treats $(u(x),p(x))$ as a
{\sl distributional or generalized} solution to
\rf{eq}--\rf{ini} in the whole
$\bbfR^3$, they
correspond to the very singular external force
$F=(b\delta_0, 0,0),$
where the parameter $b\neq 0$ depends on $c$ and $\delta_0$ stands for the
Dirac delta.

Details of this reasoning and relevant references are gathered in 
\cite{CK1}. Here, we would like only to emphasize that small 
solutions of the form \rf{sing-sol} can be also obtained from Theorem 
\ref{th:eq} because, as it is shown in Lemma \ref{F-2} below, 
$F=(c_1\delta_0, c_2\delta_0, c_3\delta_0)$ belongs to the space 
$\Y_3$ defined in \rf{Y3}.

\bigskip

The main goal of this paper is to
compare, for large $t$, properties of solutions of \rf{eq}--\rf{ini} with
properties of solutions of the following Cauchy problems:
for the  mollified Navier--Stokes system
\begin{eqnarray}
&&v_t-\Delta v +\P\nabla \cdot  ((v*\omega )\otimes v)=\P G,\label{eq-m}\\
&&v(0)=v_0.\label{ini-m}
\end{eqnarray}
where $\omega$ is a nonnegative smooth compactly supported function
on $\bbfR^3$
such that $\int_{\bbfR^3} \omega(x)\;dx =1$; and for the Navier--Stokes
system with the hyperdissipative term
\begin{eqnarray}
&&w_t-\Delta w +(-\Delta)^{\ell/2}w +\P\nabla \cdot  (w\otimes w)=\P H,
\label{eq-h}\\
&&w(0)=w_0\label{ini-h}
\end{eqnarray}
with fixed $\ell >2$.

\begin{remark}
Note that the constant $\kappa$, mentioned in Introduction,  does not 
appear in both models. In fact, without loss of generality and for 
simplicity of
notation, we put $\kappa=1$.
\end{remark}

It is not surprising that  the theories
on the existence of
global-in-time  small solutions
to  all models, \rf{eq}--\rf{ini}, \rf{eq-m}--\rf{ini-m}, and
\rf{eq-h}--\rf{ini-h} are completely analogous. Below, in Theorems
\ref{th:eq-m} and \ref{th:eq-h} , we state this fact more
precisely. However,
the main result 
of this paper consists in showing
that
the mollification of the nonlinearity in model \rf{eq-m}--\rf{ini-m}
as well as the higher order term $(-\Delta)^{\ell/2}$ with $\ell >2$
in \rf{eq-h}--\rf{ini-h} are
asymptotically negligible for large $t$. Details are contained in Theorem
\ref{asymp-m} and \ref{asymp-h}, below.

First, however, let us recall that if $u_0\in L^{3,\infty}$ is a homogeneous
function of degree $-1$ and if
$F$ satisfies
\begin{equation}
F(x,t)=\lambda^3F(\lambda x, \lambda^2 t) \quad \mbox{for all} \quad
\lambda>0,
\label{F-scal}
\end{equation}
  we already know (cf.~e.g.~\cite{CP})
that the
solution to the Navier-Stokes system \rf{eq}-\rf{ini} is
self-similar, hence of the form
\begin{equation}
U(x,t)=t^{-1/2}U(x/t^{1/2},1).\label{ss}
\end{equation}
Obviously, this is not the case
of the mollified system \rf{eq-m}--\rf{ini-m}, because it is not invariant
under the well-known rescaling $u(x,t)\mapsto \lambda
u(\lambda x,\lambda^2t)$. The goal of our next theorem is to show,
however,
that, as $t\to\infty$, solutions of \rf{eq-m}--\rf{ini-m} converge toward
suitable self-similar solutions  of the Navier-Stokes system
\rf{eq}--\rf{ini}.

\begin{theorem}\label{asymp-m}
Denote by $u=u(x,t)$ and $v=v(x,t)$ the solutions to the problems
\rf{eq}--\rf{ini} and \rf{eq-m}--\rf{ini-m}, respectively, corresponding to
the same initial datum $u_0\in L^{3,\infty}$ and the external force 
$F\in \Y_3$. Assume that
$\|u_0\|_{3,\infty}+\|F\|_{\Y_3}<\varepsilon$, $u_0$ is homogeneous of degree
$-1$, and $F$ satisfies \rf{F-scal}.
Then for  $p\in (3,\infty)$
\begin{equation}
\lim_{t\to\infty} t^{(1-3/p)/2}\|u(\cdot,t)-v(\cdot,t)\|_p=0.
\label{as-m:0}
\end{equation}
\end{theorem}

To understand the limit relation \rf{as-m:0}, one should remember 
that the self-similar solution $u=u(x,t)=t^{-1/2}U(x/\sqrt t)$ used 
in Theorem \ref{asymp-m} satisfies\newline
$
t^{(1-3/p)/2}\|u(\cdot,t)\|_p=\|U\|_p=\mbox{const.}
$
for all $t>0$ and each $p>3$. One can also look at \rf{as-m:0} in the 
following way. Let us consider the rescaled function 
$u_\lambda(x,t)=\lambda u(\lambda x, \lambda^2t)$ for each 
$\lambda>0$. Note that the self-similar solution $u$ is invariant 
under this rescaling. Hence, by a simple change of variables, for 
every fixed $t_0>0$, we obtain
\begin{eqnarray*}
  \|u(\cdot,t_0)-v_\lambda (\cdot,t_0)\|_p&=&
\|u_\lambda(\cdot,t_0)-v_\lambda (\cdot,t_0)\|_p\\
&=&\lambda^{1-3/p}\|u(\cdot,\lambda^2t_0)-v (\cdot,\lambda^2 t_0)\|_p\\
&=&(t/t_0)^{(1-3/p)/2}\|u(\cdot,t)-v (\cdot,t)\|_p
\end{eqnarray*}
after substituting $\lambda=\sqrt{t/t_0}$, $t>0$. Hence, due to these 
calculations and relation \rf{as-m:0}, it follows that
under the assumptions of Theorem \ref{asymp-m},
$v_\lambda(\cdot,t_0)\to u(\cdot,t_0)$
as $\lambda\to\infty$
in $L^p$, $3<p<\infty$, for each fixed $t_0>0$.

\medskip
The result on the asymptotic stability of the Navier-Stokes system with
hyperdissipativity is more general and reads as follows.

\begin{theorem}\label{asymp-h}
Denote by $u(x,t)$ and $w(x,t)$ the solutions to the problems
\rf{eq}--\rf{ini} and \rf{eq-h}--\rf{ini-h}, respectively, corresponding to
the same initial datum $u_0\in L^{3,\infty}$ and the  force 
of the form $F=\nabla\cdot V$ with $V\in C_w([0,\infty), L^{3/2,\infty})$.
Assume that $\|u_0\|_{3,\infty} +\|F\|_{\Y_3}\leq \varepsilon <1/(4\eta
C_\ell)$ for the constant $\eta$ defined in Proposition \ref{prop:L-S} and $C_\ell$ given by equation \rf{Cl}.
Then
$$
\lim_{t\to\infty} \|u(\cdot,t)-w(\cdot,t)\|_{3,\infty} =0.
$$
Moreover, for every $p\in (3, \infty)$,
\begin{equation}
\lim_{t\to\infty} t^{(1-3/p)/2}\|u(\cdot,t)-w(\cdot,t)\|_p=0.
\label{Lp-h}
\end{equation}
\end{theorem}

Note that $u=u(x,t)$ in Theorem \ref{asymp-h} is not assumed to be self-similar.

Section 5 contains the detailed analysis (including the proof of 
Theorem \ref{asymp-m})
of the mollified system \rf{eq-m}--\rf{ini-m}. Analogous results on the hyperviscous 
problem
\rf{eq-h}--\rf{ini-h} are gathered in Section 6. Finally, in Section 
7, we describe how
to reformulate our results in a framework of abstract functional 
Banach spaces more general than $L^p$ and $L^{p,\infty}$.

A preliminary version of  results from this paper was announced
without proof in \cite{CK0}.
The Marcinkiewicz spaces appear as well in the study of the 
Navier-Stokes system in an exterior
domain. Indeed, ``physically reasonable'' stationary solutions 
constructed by Finn \cite{Finn} in the 3-dimensional exterior
problem have the infinite energy and decay like $|x|^{-1}$ as 
$|x|\to\infty$, hence, the Marcinkiewicz space $L^{3,\infty}$ seems to be a natural space containing functions with such a behavior at infinity.
 This idea motivated to study the exterior
problem for the incompressible Navier-Stokes system in the  space 
$L^{3,\infty}(\Omega)$ (see e.g.
\cite{BM, KY-Ann,KY-Zeit,S2,Y00} and the references given there). 
Results in this direction were also
obtained in the recent paper \cite{BCK}, where ideas from 
Theorem~\ref{asymp} were adapted.


\section{Marcinkiewicz spaces}
\setcounter{equation}{0}
\setcounter{remark}{0}

In this paper, we  work in
the weak Marcinkiewicz $L^p$- spaces ($1<p<\infty$)  denoted as
usual by $L^{p,\infty}=L^{p,\infty}(\bbfR^3)$. They belong to the scale of
the
Lorentz spaces and contain measurable functions $f=f(x)$ satisfying
the condition
\begin{equation}
\left|\{x\in \bbfR^3\;:\; |f(x)|>\lambda \}\right|\leq C\lambda^{-p}
\label{18.3}
\end{equation}
for all $\lambda >0$ and a constant $C$. One can  check that
\rf{18.3} is equivalent to
$$
\int_{E} |f(x)|\;dx\leq C' |E|^{1/q}
$$
for every measurable set $E$ with a finite measure, another constant
$C'$, and $1/p+1/q=1$. This fact allows us to define the norm in
$L^{p,\infty}$
$$
\|f\|_{p,\infty} =\sup \{|E|^{-1/q} \int_E |f(x)|\;dx \;:\; E\in
{\cal B}\}
$$
where $\cal B$ is the collection of all Borel sets with a finite and
positive measure.

Recall the well-known imbedding $L^p\subset L^{p,\infty}$ being the
consequence of the Markov inequality
$|\{x\in\bbfR^n\;:\: |f(x)|>\lambda\}|\leq \lambda^{-p} \int_{\bbfR^n}
|f(x)|^p\;dx$. Moreover, in the Mar\-cin\-kie\-wicz spaces, the following
inequalities hold true: {\it the weak H\"older inequality}:
\begin{equation}
\|fg\|_{r,\infty}\leq \|f\|_{p,\infty}\|g\|_{q,\infty} \label{w-Holder}
\end{equation}
for every $1<p\leq \infty$ (here, $L^{\infty,\infty}=L^\infty$),
$1<q<\infty$, and $1<r<\infty$ satisfying $1/r=1/p+1/q$, and  {\it  the weak Young inequality}
\begin{equation}
\|f*g\|_{r,\infty}\leq C\|f\|_{p,\infty}\|g\|_{q,\infty} \label{w-Young}
\end{equation}
for every $1<p< \infty$, $1<q<\infty$, and $1<r<\infty$ satisfying $1+1/r=1/p+1/q$.

The classical Young inequality applied to the heat semigroup implies
the existence of a
constat $C=C(p,q)$ such that for every $u_0\in L^q$
\begin{equation}
\|S(t)u_0\|_p\leq Ct^{-(3/2)(1/q-1/p)}\|u_0\|_q\label{Spq},
\end{equation}
provided $1\leq q\leq p\leq \infty$. The counterpart for the Marcinkiewicz
spaces is also valid
\begin{equation}
\|S(t)u_0\|_{p,\infty}\leq
Ct^{-(3/2)(1/q-1/p)}\|u_0\|_{q,\infty}\label{Spq-w}
\end{equation}
under the additional assumption $q\neq 1$.

In the following, we also use   estimates involving the weak $L^p$ spaces
which were recently obtained independently by Y. Meyer in \cite{M96}
and by M. Yamazaki in \cite{Y00}. For the
completeness of the exposition, we  recall
them in a form most suitable for our applications.

Denote by $K(x,y,t)$, $T>0$, $x,y\in\bbfR^3$, the kernel fulfilling
the following estimate
$$
|K(x,y,t)|\leq Ct^{-3} (1+|x-y|/t)^{-4},
$$
and for every $t>0$ define the operator
$$
[P(t)h](x)=\int_{\bbfR^3} K(x,y,t) h(y)\;dy.
$$

The main estimate is contained in the following lemma.

\begin{lemma}
\label{lem:Mey} {\rm (\cite[Ch.~8, Th.~9]{M96}, \cite[Th. 3.1]{Y00})}
There exists a constant $C$ such that for every $f\in C_w([0,\infty),
L^{3/2,\infty} (\bbfR^3))$ and $g=g(x)$ defined by the formula
$$
g(x)=\int_0^\infty [P(\tau)f](x,\tau)\;d\tau
$$
we have
$$
\|g\|_{3,\infty} \leq C \sup_{\tau>0} \|f(\cdot,
\tau)\|_{3/2,\infty}.
$$\cbdu
\end{lemma}

In his consideration \cite{M96}, Meyer applied this lemma to a very
special function
\begin{equation}
g(\cdot,t)=\int_0^t \P\nabla S(t-\tau) f(\cdot,\tau)\;d\tau
\label{g}
\end{equation}
where $\P$ is the Leray projection and $S(t)$ is the heat semigroup. It is
well-known that $\P\nabla S(t)$ is given as a convolution operator with
the Oseen kernel ${\cal K}(x,y,t)=t^{-2} {\cal K}((x-y)/\sqrt t)$ where
$|{\cal K}(x)|\leq C(1+|x|)^{-4}$. Now, we change the variables $s=t-\tau$
in the
integral \rf{g} which leads to
\begin{equation}
g(\cdot,t) =\int_0^\infty \P\nabla S(s)Q(\cdot,s)\;ds
\label{g-F}
\end{equation}
with $Q(x,s) =f(x,t-s)$ if $0\leq s\leq t$, and $Q(x,s)=0$ if $s>t$.
Finally, the application of Lemma \ref{lem:Mey} gives an inequality which
plays a crucial role in our reasoning below:
\begin{equation}
\|g(\cdot,t)\|_{3,\infty} \leq \eta \sup_{0\leq s}
\|Q(\cdot,s)\|_{3/2,\infty}= \eta\sup_{0\leq \tau\leq
t}\|f(\tau,s)\|_{3/2,\infty},
\label{g:est}
\end{equation}
where the constant $\eta$ is independent of $f$, $t$, and $s$.

\medskip

Here, we also recall a result on the continuity with respect to $t$ of
$g=g(x,t)$ defined in \rf{g}.

\begin{lemma}\label{lem:Mey:2}
{\rm \cite[Lem.~24]{M96},  \cite[Th.~3.1]{Y00})}
For every $f\in C_w([0,\infty), L^{3/2,\infty}(\bbfR^3))$, the function
$g=g(x,t)$ defined in \rf{g} satisfies
$g\in C_w([0,\infty), L^{3,\infty}(\bbfR^3))$.
\cbdu
\end{lemma}

An improvement of the Meyer-Yamazaki inequality can be found in the 
recent paper by Terraneo \cite[Prop 1.5]{Ter}.

\medskip

In Section 2 (cf. \rf{Y3}), we have already defined the space $\Y_3$ 
of admissible external forces. Here, we would like to present two 
sufficient conditions for $F$ to belong to $\Y_3$.

\begin{lemma}\label{F-1}
Assume that $F(x,t)=\nabla \cdot V(x,t)$ for the external potential satisfying
$v\in C_w([0,\infty), L^{3/2,\infty})$. Then $F\in \Y_3$; moreover, if
$\lim_{t\to\infty} \|V(t)\|_{3/2,\infty}=0$, then
\begin{eqnarray}
\lim_{t\to\infty} \left\|\int_0^t S(t-\tau)\P F 
(\tau)\;d\tau\right\|_{3,\infty} =0.\label{lem:F-1:dec}
\end{eqnarray}
\end{lemma}

\proof
The first part of this Lemma is a direct consequence of Lemmata 
\ref{lem:Mey} and \ref{lem:Mey:2}.
Assume now that $\lim_{t\to\infty} \|V(t)\|_{3/2,\infty}=0$.
To prove \rf{lem:F-1:dec}, it suffices to repeat the reasoning either 
from the proof of Theorem \ref{asymp} (cf.
equation \rf{as-2}--\rf{as-3}, below) with $A={\rm 
lim\;sup}_{t\to\infty} \|V(t)\|_{3/2,\infty}=0$
or from the proof of Lemma \ref{Sl-f}. Let us skip other details.
\cbdu

\medskip

The next lemma deals with forces independent of time.

\begin{lemma}\label{F-2}
Assume that $F(\cdot,t)=\mu$ where $\mu$ is the Borel measure on 
$\bbfR^3$. Then $F\in \Y_3$.
\end{lemma}

\proof
Recall that by the definition of the space $\Y_3$, we should find an 
estimate of the norm $\|\int_0^t \P S(t-\tau)\mu 
\;d\tau\|_{3,\infty}$ which are uniform with respect to $t$. Note 
that the Leray projector $\P$ (being the combination of the Riesz 
transforms) is bounded $L^{3,\infty}$, hence it suffices to study
  $\|\int_0^t S(t-\tau)\mu \;d\tau\|_{3,\infty}$. Computing the 
Fourier transform
  of the integral $\int_0^t S(t-\tau)\mu \;d\tau$ we obtain the product
$$
\int_0^t e^{-(t-\tau)|\xi|^2}\;d\tau \;\widehat \mu(\xi)
= {1-e^{-t|\xi|^2}\over |\xi|^2}\widehat \mu(\xi).
$$
Hence, $\int_0^t S(t-\tau)\mu \;d\tau=E_3*((I-S(t))\mu)$, where
$E_3(x)=(4\pi |x|)^{-1}$  is the fundamental 
solution of the Laplace operator on $\bbfR^3$. Since $E_3\in
L^{3,\infty}$ and convolutions of Borel measures with elements from 
$L^{3,\infty}$ are well-defined, we obtain
$$
\|E_3*((I-S(t))\mu)\|_{3,\infty}\leq C\|E_3\|_{3,\infty}<\infty.
$$
We skip the proof of the regularity with respect to $t$ because the 
reasoning is more or less similar to that used in the proof of Lemma 
\ref{lem:Mey:2}.
\cbdu
\section{The Navier-Stokes system}
\setcounter{equation}{0}
\setcounter{remark}{0}

As in \cite{C}, the proof of our  theorem
on the existence, uniqueness and stability of solutions to the
problem \rf{eq}--\rf{ini} is based on the following abstract lemma,
whose slightly more general form is taken
from \cite{LR}.

\begin{lemma}
\label{lem:xyB}
Let $(\X, \|\cdot\|_\X)$ be a Banach space and $B:\X\times
\X\to
\X$ a bounded bilinear form satisfying
$
\|B(x_1,x_2)\|_\X\leq \eta \|x_1\|_\X \|x_2\|_\X
$ for all $x_1,x_2\in\X$ and a constant $\eta>0$.
Then, if $0<\varepsilon<1/(4\eta)$ and if $y\in\X$  such that
$ \|y\|<\varepsilon$, the equation
$x=y+B(x,x)$ has a solution in $\X$ such that $\|x\|_\X\leq
2\varepsilon$.
This solution is the
only one in the ball $\bar B(0,2\varepsilon)$.
Moreover, the solution depends continuously on $y$ in the following
sense: if $\|\tilde y\|_\X\leq \varepsilon$,
$\tilde x=\tilde y+B(\tilde x, \tilde x)$, and $\|\tilde x\|_\X\leq
2\varepsilon$, then
$$
\|x-\tilde x\|_\X \leq {1\over 1-4\eta \varepsilon}\|y-\tilde y\|_\X.
$$
\end{lemma}

\proof Here, one uses the standard Picard iteration
technique completed by the Banach fixed point theorem. For other
details of the proof, we refer the reader to \cite[Th. 13.2]{LR}.
\cbdu
\medskip

Our goal is to apply Lemma \ref{lem:xyB}
to the integral equation \rf{duhamel}
in the space $\X_3$ defined in \rf{X3}.
To continue, we need the
estimate of the  form $B(\cdot,\cdot)$.

\begin{proposition}\label{prop:L-S}
The bilinear form $B(\cdot,\cdot)$ is bounded on the space
$\X_3$. In other words, there
exists a constant $\eta>0$ such that
$$
\|B(u,v)\|_{\X_3}\leq \eta \|u\|_{\X_3} \|v\|_{\X_3}
$$
for all $u,v\in \X$.
\end{proposition}

\proof The proof of this fact,  given by Meyer in
\cite[Ch.~18]{M96},  results immediately from Lemmata \ref{lem:Mey}
and \ref{lem:Mey:2}
because, for all $u,v\in \X$, it follows that $u\otimes v \in
C_w([0, \infty), L^{3/2,\infty})$. An independent reasoning which
leads to this proposition can be also found in the recent paper by M.
Yamazaki \cite[Th.~3.1]{Y00}.
\cbdu

\medskip

\noindent {\sc Proof of Theorem \ref{th:eq}.}
Now, the main theorem on the existence of unique small solutions  is 
a consequence of
Lemma \ref{lem:xyB} combined with  Proposition \ref{prop:L-S}.\cbdu
\medskip

\begin{remark}
Homogeneity properties of equation \rf{eq} imply that if
$u$ solves the Cauchy problem, then the rescaled function
$u_\lambda(x,t)=\lambda u(\lambda x,\lambda^2t)$ is also a solution
for each $\lambda>0$. Thus, it is natural to consider solutions which
satisfy the scaling invariance property $u_\lambda\equiv u$ for all
$\lambda>0$, i.e. {\sl forward self-similar } solutions.
By the uniqueness property of solutions of
the Cauchy problem,
they can be obtained directly from Theorem
\ref{th:glob}
by taking $u_0$ homogeneous of degree $-1$  and $F$ satisfying
\rf{F-scal}.
\end{remark}

Regularity of  solutions constructed in Theorem \ref{th:eq} depends
essentially on the regularity of external forces. We precise this fact in
our next proposition.

\begin{proposition}\label{eq:reg}
Let the assumptions of Theorem \ref{th:eq} hold true. 
Assume that  $3<p\leq \infty$ and recall that the space $\X_p$ is defined in \rf{Xp}.
Suppose, moreover,
that the external force $F$ satisfies $\int_0^t S(t-\tau)\P F(\tau)\;d\tau\in \X_p$. 
There exists $\varepsilon_p\in
(0,\varepsilon]$ such that the solution $u=u(x,t)$ constructed in Theorem
\ref{th:eq}
belong to the space $\X_3\cap \X_p$. 
\cbdu
\end{proposition}

We skip the proof of this proposition, because it is more or less standard
(see e.g. \cite{Ba1,Ba2,C,Ca2,CPfo,Kat1,K,LR,Y00}, for details). Let us only 
mention that it is based on Lemma~\ref{lem:xyB} applied in the space $\X=\X_3\cap\X_p$, and the 
required estimates
of the bilinear form $B(\cdot,\cdot)$ defined in \rf{Bf} can be easily
obtained combining  the well-known inequalities for the heat 
semigroup \rf{Spq} and its derivatives with the H\"older inequality. 
Here, the crucial role is played by the inequality
\begin{equation}
\|\nabla\P S(t-\tau)(f\otimes g)\|_p\leq \eta_p
(t-\tau)^{-(1+3/p)/2}\|f\|_p\|g\|_p \label{Sfgp}
\end{equation}
valid
for every $p\in (3,\infty]$
all $0<\tau<t$ and a constant $\eta_p$, as well as its counterpart in the $L^{p,\infty}$-spaces. Note that \rf{Sfgp} holds also true for
$p\in [2,3]$ but, in this case, the function
$\zeta(\tau)=(t-\tau)^{-(1+3/p)/2}$ is not integrable near $\tau=t$.
All details concerning the proof of Proposition \ref{eq:reg} are 
contained in \cite[Th.~5.1]{K}, \cite[Th. 1]{Ba1}.

\bigskip

\noindent {\sc Proof of Theorem \ref{asymp}.}
Several estimates from this proof will be used later on in the
  analysis of the regularized problems \rf{eq-m}--\rf{ini-m} and
  \rf{eq-h}--\rf{ini-h}, hence we shall try to be very detailed.

We begin by recalling that, by Theorem \ref{th:glob}, we have
\begin{equation}
\sup_{t\geq 0} \|u(t)\|_{3,\infty} \leq 2\varepsilon<{1\over 2\eta}
\quad
\mbox{and} \quad
     \sup_{t\geq 0} \|\tilde u(t)\|_{3,\infty}\leq
     2\varepsilon<{1\over 2\eta}.
\label{5.1.u.v}
\end{equation}
We subtract the integral equation \rf{duhamel} for $\tilde u$ from
the analogous expression for $u$. Next, computing the norm
$\|\cdot\|_{3,\infty}$ of the resulting equation  we obtain the following
inequality
\begin{eqnarray}
&&\hspace{-1cm}\|u(t)-\tilde u(t)\|_{3,\infty} \nonumber\\
&\leq& \left\|S(t)(u_0-\tilde u_0)
+\int_0^t S(t-\tau)\P (F-\tilde
F)(\tau)\;d\tau
\right\|_{3, \infty} \label{as-1}\\
&&  +\left\| \int_0^{\delta t} S(t-\tau) \P\nabla \cdot
[(u-\tilde u)\otimes u+\tilde u \otimes (u-\tilde u)](\tau)\;d\tau
\right\|_{3,\infty}
\nonumber\\
&& +\left\| \int_{\delta t}^t S(t-\tau) \P\nabla \cdot
[(u-\tilde u)\otimes u+\tilde u \otimes (u-\tilde u)](\tau)\;d\tau
\right\|_{3,\infty}\nonumber
\end{eqnarray}
where the small constant $\delta>0$ will be chosen later.

In the term on the right-hand side  of \rf{as-1}
containing the integral $\int_0^{\delta t}
...\;d\tau$, we apply the weak $L^p-L^q$ estimates of the heat semigroup
\rf{Spq-w}, the boundedness of $\P$ on $L^{3,\infty}$, the weak 
H\"older inequality \rf{w-Holder}, and \rf{5.1.u.v},
in order to estimate it by
\begin{eqnarray}
&&  C  \int_0^{\delta t}  (t-\tau)^{-1}\|u(\tau)-\tilde
u(\tau)\|_{3,\infty}\;
d\tau \nonumber\\
&&\hspace{1.5cm}\times
\left(\sup_{\tau>0} \|u(\tau)\|_{3,\infty}+\sup_{\tau>0}
\|\tilde u(\tau)\|_{3,\infty}
\right)\label{as-2}\\
&& \hspace{1cm}
\leq 4\varepsilon C  \int_0^{\delta} (1-s)^{-1}\nonumber
\|u(ts)-v(ts)\|_{3,\infty}\; ds.
\end{eqnarray}

To deal with the term in \rf{as-1}
     containing  $\int_{\delta t}^t ...\;d\tau$, we use Lemma \ref{lem:Mey}
(with $f=(u-\tilde u)*u+\tilde u*(u-\tilde u)$ for $\delta t<\tau< t$ 
and $f=0$ otherwise) combined with the H\"older inequality 
\rf{w-Holder} and with \rf{5.1.u.v},  to
  bound it directly by
\begin{eqnarray}
&& \eta \left(
\sup_{\delta t\leq \tau\leq t} \|u(\tau)\|_{3,\infty}
+
\sup_{\delta t\leq \tau\leq t} \|\tilde u(\tau)\|_{3,\infty}
\right)
     \sup_{\delta t\leq \tau\leq t} \|u(\tau) -\tilde u
     (\tau)\|_{3,\infty}
\label{as-3}\\
&&\hspace{1cm}
\leq  4\varepsilon \eta \sup_{\delta t\leq \tau\leq t} \|u(\tau)
-\tilde u(\tau)\|_{3,\infty}.
\nonumber
\end{eqnarray}

Now, we denote
$$
g(t)= \left\|S(t)(u_0-\tilde u_0)
+\int_0^t S(t-\tau)\P (F-\tilde
F)(\tau)\;d\tau
\right\|_{3,\infty},
$$
and it follows from   the assumptions on initial data and external forces
  that
\begin{equation}
g\in L^\infty(0, \infty) \quad \mbox{and} \quad
\lim_{t\to\infty} g(t)=0. \label{as-3a}
\end{equation}
Hence, applying \rf{as-2} and \rf{as-3} to \rf{as-1} we arrive at
\begin{eqnarray}
\|u(t)-\tilde u(t)\|_{3,\infty} &\leq & g(t) +
4\varepsilon C  \int_0^{\delta} (1-s)^{-1}
\|u(ts)-\tilde u (ts)\|_{3,\infty}\; ds \label{as-4}\\
&& + 4\varepsilon \eta \sup_{\delta t\leq \tau\leq t} \|u(\tau)
-\tilde u(\tau)\|_{3,\infty}
\nonumber
\end{eqnarray}
for all $t>0$.

Next, we put
$$
A=\limsup_{t\to\infty} \|u(t)-\tilde u(t)\|_{3,\infty} \equiv
\lim_{k\in \bbfN, k\to\infty} \sup_{t\geq k} \|u(t)-\tilde u
(t)\|_{3,\infty}.
$$
The number $A$ is nonnegative and finite because both
$u,\tilde u\in  L^\infty ([0,\infty), L^{3,\infty})$, and our claim
is to show that
$A=0$.

First, we apply  the Lebesgue dominated convergence theorem to the
obvious inequality
$$
\sup_{t\geq k} \int_0^{\delta} (1-s)^{-1} \|u(ts)-\tilde u
(ts)\|_{3,\infty} \;ds
\leq \int_0^{\delta} (1-s)^{-1} \sup_{t\geq k} \|u(ts)-\tilde
u(ts)\|_{3,\infty}
\;ds,
$$
and we obtain
\begin{equation}
\limsup_{t\to\infty}\int_0^{\delta} (1-s)^{-1} \|u(ts)-\tilde
u(ts)\|_{3,\infty}
\;ds
\leq A\int_0^{\delta} (1-s)^{-1}\;ds =A\log \left({1\over
1-\delta}\right). \label{as-5}
\end{equation}
Moreover, since
$$
\sup_{t\geq k} \sup_{\delta t\leq\tau\leq t} \|u(\tau)
-\tilde u(\tau)\|_{3,\infty}
\leq \sup_{\delta k\leq \tau <\infty} \|u(\tau)-\tilde
u(\tau)\|_{3,\infty},
$$
we have
\begin{equation}
\limsup_{t\to\infty}  \left( \sup_{\delta t\leq \tau\leq t} \|u(\tau)
-v(\tau)\|_{3,\infty}\right)\leq A.\label{as-6}
\end{equation}

Finally, computing $\limsup_{t\to\infty}$ of the both sides of
inequality \rf{as-4}, and using \rf{as-3a}, \rf{as-5}, and \rf{as-6}
we get
$$
A\leq \left( 4 \varepsilon C\log \left({1\over
1-\delta}\right) +4
\varepsilon\eta\right)A.
$$
Consequently,
it follows that
$
A=\limsup_{t\to\infty}\|u(t)-v(t)\|_{3,\infty}=0
$
because
$$
4 \varepsilon \eta \left(C\log \left({1\over 1-\delta}\right)
+1\right) <1,
$$
for $\delta>0$ sufficiently small,
by the assumption of Theorem \ref{th:glob} saying that $0<
\varepsilon <1/(4\eta)$.
This completes the proof of Theorem \ref{asymp}.
\cbdu

\medskip

As a direct consequence the proof of Theorem \ref{asymp},
we have also necessary conditions for \rf{e:utvt} to hold. We
formulate this
fact
in the following corollary.

\begin{corollary}\label{or:as}
Assume that $u,\tilde u \in \X_3$ are solutions to system
\rf{eq}--\rf{ini}
corresponding to initial conditions $u_0, \tilde u_0\in
L^{3,\infty}$ and external forces $F,\tilde F\in \Y_3$, respectively.
Suppose that
\begin{equation}
\lim_{t\to\infty} \|u(t)-\tilde u(t)\|_{3,\infty}=0. \label{c0}
\end{equation}
Then
$$
\lim_{t\to\infty} \left\|S(t)(u_0-\tilde u_0)
+\int_0^t S(t-\tau)\P
(F(\tau)-\tilde F(\tau))\;d\tau\right\|_{3,\infty}=0.
$$
\end{corollary}

\proof
As in the beginning of the proof of Theorem \ref{asymp}, we subtract
the integral equation
\rf{duhamel} for $\tilde u$ from the same expression for $u$. Next, we
compute
the $L^{3,\infty}$-norm and we use inequalities \rf{as-2} and \rf{as-3} to
obtain
\begin{eqnarray}
&& \hspace{-1cm}\left\|S(t)(u_0-\tilde u_0)+\int_0^t S(t-\tau)\P
(F(\tau)-\tilde F(\tau))\;d\tau\right\|_{3,\infty}
\nonumber\\
&& \leq\|u(t)-\tilde u(t)\|_{3,\infty}\label{c1}\\
&&\hspace{0.5cm}+
4\varepsilon C  \int_0^{\delta} (1-s)^{-1}
\|u(ts)-v(ts)\|_{3,\infty}\; ds
+4\varepsilon \eta \sup_{\delta t\leq \tau\leq t} \|u(\tau)
-\tilde u(\tau)\|_{3,\infty}.
\nonumber
\end{eqnarray}
The first term on the right-hand side of \rf{c1} tends to zero as
$t\to\infty$
by  \rf{c0}. To show the decay of  the second one, it suffices to
repeat calculations
from  \rf{as-3} and \rf{as-5}.  Now, however,
one should remember that $A=0$ is assumed. \cbdu

\bigskip

An asymptotic stability result holds also true in the $L^{p,\infty}$ and
$L^p$-spaces with $p>3$.

\begin{proposition}\label{prop:Lp}
Under the assumptions of Theorem \ref{asymp}
$$
\lim_{t\to\infty} t^{(1-3/p)/2}\|u(\cdot,t)-\tilde
u(\cdot,t)\|_{p,\infty}=0
$$
and
$$
\lim_{t\to\infty} t^{(1-3/q)/2}\|u(\cdot,t)-\tilde
u(\cdot,t)\|_{q}=0
$$
for every $q\in (3,p)$.
\cbdu
\end{proposition}

We skip the proof of the first relations in the above proposition, 
because it a standard fact. A completely analogous reasoning can be 
found  in
  \cite{P2, K, Ba2} (see also the last section of this paper).
Next, the limit in the $L^q$-spaces is a consequence of the imbedding 
\rf{imbedding}.

\section{The mollified Navier-Stokes system}
\setcounter{equation}{0}
\setcounter{remark}{0}

The goal of this section is to formulate and to prove results on the
large time behavior of solutions to the mollified problem
\rf{eq-m}--\rf{ini-m} reformulated in the integral form
\begin{eqnarray}
v(t)= S(t) v_0 +B_\omega(v,v)
  +\; \int_0^t S(t-\tau) \P G(\tau)\;d\tau \label{duhamel-m}
\end{eqnarray}
with the bilinear form
$$
B_\omega (v, \tilde v)(t)=-
\int_0^t S(t-\tau) \P\nabla \cdot [(v*\omega)\otimes \tilde
v](\tau)\;d\tau.
$$

The counterpart of Proposition \ref{prop:L-S} reads as follows.

\begin{proposition}
For every $v, \tilde v\in \X_3$, we
have $B_\omega (v,\tilde v)\in \X_3$.
Moreover, it follows that
$$
\|B_\omega (v,\tilde v)\|_{\X_3}\leq \eta \|v\|_{\X_3}\|\tilde v\|_{\X_3}
$$
with the same constant  $\eta$   as in Proposition \ref{prop:L-S}.
\end{proposition}

\proof
This is a direct application of Proposition \ref{prop:L-S}, since by the
H\"older inequality, we have
$
\|v*\omega \|_{\X_3}\leq \|v\|_{\X_3} \|\omega \|_{L^1} = \|v\|_{\X_3}.
$
\cbdu

\medskip

Now, as usual, the existence, uniqueness, and regularity  of solutions to
\rf{eq-m}--\rf{ini-m} are  deduced from Lemma \ref{lem:xyB}.

\begin{theorem}\label{th:eq-m}
The counterpart of Theorem \ref{th:glob} for the mollified problem
\rf{eq-m}--\rf{ini-m}
holds true if we replace  $u(x,t)$,
$u_0$, and $F$ by the solution $v=v(x,t)$ to
\rf{eq-m}--\rf{ini-m} corresponding to the initial datum $v_0\in L^{3,\infty}$ and the
external force $G\in\Y_3$.
\cbdu
\end{theorem}

\begin{proposition}\label{eq:reg-m}
Let the assumptions of Theorem \ref{th:eq-m} hold true. 
Let $3<p\leq \infty$.
Suppose, moreover,
that the external force $G$ satisfies $\int_0^t S(t-\tau)\P G(\tau)\;d\tau\in \X_p$. 
There exists $\varepsilon_p\in
(0,\varepsilon]$ such that the solution $v=v(x,t)$ constructed in Theorem
\ref{th:eq-m}
belong to the space $\X_3\cap \X_p$. 
\cbdu
\end{proposition}

Here, we have skiped the proofs of  Theorem \ref{th:eq-m} and Proposition \ref{eq:reg-m} because they are completely analogous to their counterparts 
from Section 4.

Now, we are in a position to prove the convergence of solutions of 
the mollified problem \rf{eq-m}--\rf{ini-m} toward
self-similar solutions of \rf{eq}--\rf{ini}.

\medskip

\noindent {\sc Proof of Theorem \ref{asymp-m}.}
Recall that by Propositions \ref{eq:reg} and  \ref{eq:reg-m}, 
$u$ and $v$ exist
for all $t>0$, they both belong to the ball $B(0,2\varepsilon_p)\subset
\X_p =\{ u\in C((0,\infty), L^p(\bbfR^3))\;:\;
\sup_{t>0} t^{(1-3/p)/2} \|u(t)\|_p<\infty\}$
for every $p\in (3,\infty]$ and some $\varepsilon_p\in (0,1/(4\eta_p)$
which implies that
\begin{equation}
\|u\|_{p}<2\varepsilon_p t^{(1-3/p)/2}\quad \mbox{and}\quad
\|v\|_{p}<2\varepsilon_pt^{_(1-3/p)/2}.\label{a-m-0}
\end{equation}
Moreover, $u=u(x,t)$ is the self-similar solution of the form \rf{ss}.

Here, we study again the difference of the integral formulations of both
problems (see \rf{duhamel} and \rf{duhamel-m}) written in the following
form
\begin{eqnarray}
u(t)-v(t)&=& -\int_0^t \nabla\P S(t-\tau) [(u-u*\omega)\otimes
u](\tau)\;d\tau\label{a-m-1}\\
&&-\int_0^t \nabla\P S(t-\tau) [((u-v)*\omega)\otimes u](\tau)\;d\tau.
\nonumber
\end{eqnarray}

It follows from the self-similar form of $u(x,t)$ that
\begin{eqnarray}
(u*\omega)(x,t)&=& t^{-1/2}\int_{\bbfR^3} U\left({x-y\over \sqrt t}\right)
\omega (y)\;dy\nonumber\\
&=& t^{-1/2}\int_{\bbfR^3} U\left({x\over \sqrt t}-z\right) t^{3/2}\omega
(y\sqrt{t})\;dy\label{a-m-3}\\
&=& t^{-1/2} (U*\omega^t)\left({x\over \sqrt t}\right),\nonumber
\end{eqnarray}
where $\omega^t(z)=t^{3/2}\omega(z\sqrt t)$. One can easily check that
$\omega^t$ is the approximation of the Dirac delta as $t\to\infty$. Hence,
in particular,
\begin{equation}
\|U-U*\omega^t\|_p\to 0\quad \mbox{as}\quad t\to\infty \label{a-m-4}
\end{equation}
for every $p\in (3,\infty)$, because $U\in L^p$ in this range of $p$.

Now, we compute the $L^p$-norm of \rf{a-m-1}, next, we multiply the
resulting inequality by $t^{(1-3/p)/2}$ and, finally, we use 
inequality \rf{Sfgp} in
order to obtain
\begin{eqnarray}
&&\hspace{-1cm}t^{(1-3/p)/2}\|u(t)-v(t)\|_p\nonumber\\
&\leq&
\eta_p t^{(1-3/p)/2}\int_0^t (t-\tau)^{-(1+3/p)/2}
\|(u-u*\omega)(\tau)\|_p \|u(\tau)\|_p\;d\tau\nonumber\\
&\leq&
\eta_p t^{(1-3/p)/2}\int_0^t (t-\tau)^{-(1+3/p)/2}
\|[(u-v)*\omega](\tau)\|_p \|v(\tau)\|_p\;d\tau\label{a-m-2}\\
&=& C(t)+D(t).\nonumber
\end{eqnarray}

The first term on the right-hand-side of \rf{a-m-2} tends to 0  as
$t\to\infty$. To see this fact, we use \rf{a-m-0}, \rf{a-m-3}, and the
change of variables $\tau=ts$ in order to show that
\begin{eqnarray*}
C(t)&\leq& 2\varepsilon_p\eta_p t^{(1-3/p)/2}\int_0^t
(t-\tau)^{-(1+3/p)/2} \tau^{-(1-3/p)}\|U-U*\omega^\tau\|_p\;d\tau\\
&=& 2\varepsilon_p\eta_p \int_0^1 (1-s)^{-(1+3/p)/2}
s^{-(1-3/p)}\|U-U*\omega^{ts}\|_p\;d\tau.
\end{eqnarray*}
Now, $\lim_{t\to\infty}C(t)=0$ by the Lebesgue dominated convergence
theorem.

We apply a similar argument involving \rf{a-m-0}
and the change of variables $\tau=ts$ to estimate the second term in
\rf{a-m-2}
by
\begin{equation}
D(t)\leq 2\varepsilon_p\eta_p \int_0^1 (1-s)^{-(1+3/p)/2}
s^{-(1-3/p)}\left((st)^{(1-3/p)/2}\|u(ts)-v(ts)\|_p\right)\;d\tau.
\label{a-m-B}
\end{equation}

Next, we define the number
$$
A=\limsup_{t\to\infty} t^{(1-3/p)/2}\|u(t)-v(t)\|_p \equiv
\lim_{k\in \bbfN, k\to\infty} \sup_{t\geq k} t^{(1-3/p)/2} \|u(t)-v(t)\|_p
$$
which is nonnegative and finite because $u,v\in X_p$, and our claim is to
show that $A=0$. Since $\limsup_{t\to\infty} C(t)=0$, it follows from
\rf{a-m-2}, \rf{a-m-B}, and from the Lebesgue dominated convergence
theorem that
\begin{equation}
A\leq \left(
2\varepsilon_p\eta_p \int_0^1 (1-s)^{-(1-3/p)/2}s^{-(1-3/p)}\;ds
\right)A
\label{a-A}
\end{equation}
The quantity in the parentheses is smaller than 1 provided  
$\varepsilon_p$ is sufficently small. Hence, inequality \rf{a-A} implies that $A=0$ and the
proof   of Theorem \ref{asymp-m} is complete.

\cbdu

\section{The hyperviscous Navier-Stokes system}
\setcounter{equation}{0}
\setcounter{remark}{0}

In the case of the system \rf{eq-h}--\rf{ini-h}, the counterpart of
the integral
equation \rf{duhamel} has the following form
\begin{eqnarray}
u(t)&=& S_\ell(t)S(t) u_0 - \int_0^t  S_\ell(t-\tau)S(t-\tau)
\nabla \P \cdot (u\otimes
u)(\tau)\;d\tau. \label{duhamel-h}\\
&&+\; \int_0^t S_\ell(t-\tau) S(t-\tau) \P H(\tau)\;d\tau, \nonumber
\end{eqnarray}
where  the semigroup generated by the operator
$(-\Delta)^{\ell/2}$ is denoted by
$S_\ell(t)$ which  is given by the convolution with the kernel
$$
p_\ell(x,t)= \int_{\bbfR^n} e^{-t|\xi|^\ell +ix\cdot \xi}\;d\xi.
$$
Note that $p_2(x,t)$ corresponds to the Gauss-Weierstrass kernel
$p(x,t)$.
  Recall  that the function
$$
p_\ell(x,t) = \int_{\bbfR^n} e^{-t|\xi|^\ell +ix\cdot \xi}\;d\xi=
t^{-n/\ell} p(x/t^{1/\ell},1)
$$
is integrable for every $\ell>0$ and all $t>0$. Moreover, the
self-similar form of $p_\ell$ implies that
$
\|p_\ell (\cdot, t)\|_{1} =
\|p_\ell (\cdot, 1)\|_{1}
$
for every $t>0$.
In this section, the constant
\begin{equation}
C_\ell=\|p_\ell (\cdot, 1)\|_{1}\label{Cl}
\end{equation}
appears
quite often in our calculation
because of the inequality
\begin{equation}
\|S_\ell(t)h\|_{3,\infty} \leq C_\ell \|h\|_{3,\infty}
\label{Cl-est}
\end{equation}
valid for every $h\in L^{3,\infty}$
(cf. also Proposition \ref{prop-hyp}, below).

\begin{remark}
Since the kernel $p_\ell(\cdot,1)$ is  integrable  for every $\ell>0$, the
constant $C_\ell$  is well-defined. However, only for $0<\ell\leq 2$, it
is a nonnegative function, consequently, in this range of $\ell$, we have
$
C_\ell=\|p_\ell(\cdot,1)\|_1=\int_{\bbfR^3} p_\ell(x,1)\;dx= \widehat
p_\ell (0,1)=1.
$
On the other hand, the kernel $p_\ell(x,t)$ changes sign for $\ell>2$,
hence, for those $\ell$, $C_\ell>1$.
\end{remark}

We  define the bilinear form
$$
B_\ell (w, \tilde w)(t)=-
\int_0^t S_\ell (t-\tau) S(t-\tau) \P\nabla \cdot [w\otimes \tilde
w](\tau)\;d\tau.
$$

The following proposition plays again an essential role in our proofs of
existence of
global-in-time solutions to \rf{eq-h}--\rf{ini-h} as well as in the study
of their large time asymptotics.

\begin{proposition}\label{prop-hyp}
For every $w, \tilde w\in \X_3$, we have $B_\ell (w,\tilde w)\in \X_3$.
Moreover, it follows that
$$
\|B_\ell (w,\tilde w)\|_{\X_3}\leq \eta C_\ell \|w\|_{\X_3}\|\tilde w\|_{\X_3}
$$
for $\eta$ defined in Proposition \ref{prop:L-S} and $C_\ell$ given by
\rf{Cl}.
\end{proposition}

\proof
This inequality results immediately from the Meyer-Yamazaki estimate
\rf{g:est} applied to the function defined in \rf{g-F} with
$Q(\cdot,\tau)=S_\ell(t-\tau) (w(\tau)\otimes\tilde w(\tau))$ for $0<\tau
<t$ and $Q(\cdot,\tau)=0$ otherwise. Next, one should use inequalities
\rf{Cl-est} and \rf{w-Holder} in the following way
\begin{eqnarray*}
\sup_{\tau>0}\|Q(\cdot,\tau)\|_{3,\infty} &=&
\sup_{0<\tau<t} \|S_\ell(t-\tau) (w(\tau)\otimes\tilde
w(\tau))\|_{3,\infty} \\
&\leq& C_\ell
\sup_{0<\tau<t}\|w(\tau)\|_{3,\infty}
\sup_{0<\tau<t}\|\tilde w(\tau)\|_{3,\infty}.
\end{eqnarray*}
\cbdu

\begin{theorem}\label{th:eq-h}
The counterpart of Theorem \ref{th:glob} holds true if we replace
$u(x,t)$, $u_0$, and $F$  by the solution $w=w(x,t)$ to the hyperviscous
problem
\rf{eq-h}--\rf{ini-h} corresponding to the initial datum $w_0\in L^{3,\infty}$ and the
external force $H\in \Y_3$, and if we
impose additional assumption
$\varepsilon <1/(4\eta C_\ell)$.
\cbdu
\end{theorem}

\begin{remark}
As in the case of problems \rf{eq}-\rf{ini} and \rf{eq-m}-\rf{ini-m}, the solution constructed in Theorem \ref{th:eq-h} belongs to the space $\X_p$ for $p>3$ under the additional assumption $\int_0^t S_\ell(t-\tau) S(t-\tau) \P H(\tau)\;d\tau \in \X_p$. Here, we omit details because the reasoning is completely analogous to that used in Propositions \ref{eq:reg} and 
\ref{eq:reg-m}.
\end{remark}

The crucial lemma in the study of the large time behavior of solutions to
\rf{eq-h}--\rf{ini-h}
  says that the
semigroup generated by the operator $\Delta -(-\Delta)^{\ell/2}$ can
be well-approximated in $L^1$ by the heat semigroup $S(t)$.

\begin{lemma} \label{SaS}
Let $\ell>0$.  There exists a constant $C$
independent of $t$ such that
$$
\|p_\ell(t)*p(t/2)-p(t/2)\|_{1} \leq Ct^{-(1/2-1/\ell)}
$$
for all $t>0$.
\end{lemma}

\proof
Let us recall
the inequality
\begin{equation}
\left\|f*g(\cdot) - \left(\int_{\bbfR^n} f(x)\;dx\right)
g(\cdot)\right\|_{1}
\leq C\|\nabla g\|_1 \|f\|_{L^1(\bbfR^n,\; |x|\,dx)}
\label{f-g-est}
\end{equation}
which is valid for all sufficiently regular $f$ and $g$, and a
constant $C$
independent of $f$, $g$. The proof of \rf{f-g-est} 
(based on the Taylor expansion of the function $g$) and its
generalizations can be found in \cite{DZ92}. Now,
in \rf{f-g-est}, we substitute
$$
f(x)=p_\ell(t,x)\quad \mbox{and} \quad
g(x,t)=p(t/2,x)
$$
to obtain (recall that $\int_{\bbfR^n} p_\ell(t,x)\;dx=\widehat
p(0,t)=1$)
\begin{eqnarray*}
\|p_\ell(t)*p(t/2)-p(t/2)\|_1
&\leq&   C\|\nabla p(t/2)\|_1
\|p_\ell(t)\|_{L^1(\bbfR^n,\;
|x|\,dx)}\\
&=& Ct^{-(1/2-1/\ell)}
\end{eqnarray*}
for all $t>0$.
\cbdu

\medskip

The lemma above is used in the proof of our next result.

\begin{lemma}\label{Sl-f}
Assume that $f\in C_w([0, \infty), L^{3/2,\infty}(\bbfR^3))$.
Then
\begin{equation}
\left\|\int_0^t \nabla\P S(t-\tau) \bigg(S_\ell(t-\tau)-I\bigg)
f(\cdot,\tau)\;d\tau\right\|_{3,\infty} \to 0 \quad \mbox{as}
\quad t\to\infty,\label{Slf-0}
\end{equation}
where $I$ denotes the identity operator.
\end{lemma}

\proof
First note that the quantity in \rf{Slf-0} is bounded uniformly with
respect to $t>0$ in view of
inequality \rf{g:est} (cf. also the proof of Proposition \ref{prop-hyp}).
To show its convergence
to 0 we fix $\gamma\in (0,1)$ (to be chosen later on) and we decompose the
integral in \rf{Slf-0} as $\int_0^t...\;d\tau =\int_0^{\gamma t}
...\;d\tau +\int_{\gamma t} ^t...\;d\tau$. Now, the estimates of the
heat semigroup on the Marcinkiewicz spaces \rf{Spq-w} give
\begin{eqnarray}
&&\hspace{-1cm}\left\|\int_0^{\gamma t} \nabla\P S(t-\tau)
\bigg(S_\ell(t-\tau)-I\bigg)
f(\cdot,\tau)\;d\tau\right\|_{3,\infty}\nonumber\\
&&\leq
\int_0^{\gamma t} \left(\left\| \nabla\P S(t-\tau) S_\ell(t-\tau)
f(\cdot,\tau)\right\|_{3,\infty}+
\left\| \nabla\P  S(t-\tau) f(\cdot,\tau)\right\|_{3,\infty}
\right)\;d\tau\nonumber\\
&& \leq C \int_0^{\gamma t} (t-\tau)^{-1}\;d\tau \;\sup_{\tau>0}
\|f(\cdot,\tau)\|_{3/2,\infty}\label{lA}\\
&&\leq C \log\left({1\over 1-\gamma}\right)\;\sup_{\tau>0}
\|f(\cdot,\tau)\|_{3/2,\infty}.\nonumber
\end{eqnarray}
Observe that the right-hand-side of the inequality above can be made
arbitrarily small choosing $\gamma>0$ sufficiently small.

We handle the integral over $[\gamma t, t]$ using the Meyer-Yamazaki
estimate \rf{g:est} applied to the function
$$
g(\cdot,t)=\int_0^t \nabla\P S\left({t-\tau\over 2}\right)Q(\cdot,
\tau)\;d\tau
$$
with $Q(\cdot,\tau)= [S_\ell(t-\tau)-I]S((t-\tau)/2) f(\cdot,\tau)$ for
$\tau \in [\gamma t,t]$ and $Q(\cdot,\tau)=0$ otherwise.
First, using  inequality \rf{g:est} and next, Lemma \ref{SaS} we obtain
\begin{eqnarray*}
\|g(\cdot,t)\|_{3,\infty} &\leq& C \sup_{\gamma t\leq \tau\leq t}
\|Q(\cdot,\tau)\|_{3/2,\infty} \\
&\leq& C \left( \sup_{\gamma t\leq \tau\leq t} (t-\tau)^{-(1/2-1/\ell)}
\right)
\sup_{\gamma t\leq \tau\leq t}\|f(\cdot,\tau)\|_{3/2,\infty}\\
&\leq& C t^{-(1/2-1/\ell)} \sup_{\gamma t\leq \tau\leq t}
\|f(\cdot,\tau)\|_{3/2,\infty}
\end{eqnarray*}
Note now that the right-hand-side of the above inequality tends to 0 as
$t\to\infty$ for every $\gamma >0$.
This completes the proof of Lemma \ref{Sl-f}.
\cbdu

\bigskip

\noindent{\sc Proof of Theorem \ref{asymp-h}.}
Note first  that the existence of such solutions is provided by Theorems
\ref{th:eq} and \ref{th:eq-h}. In particular, we have
\begin{equation}
\sup_{t\geq 0} \|u(t)\|_{3,\infty} \leq 2\varepsilon
\quad
\mbox{and} \quad
     \sup_{t\geq 0} \|\tilde w(t)\|_{3,\infty}\leq
     2\varepsilon.
\end{equation}

Here, we  describe only how to modify the proof of Theorem \ref{asymp}.
First, we subtract the integral equation \rf{duhamel-h} for $w$ from  equation
\rf{duhamel} for $u$ and,
next, we compute the $L^{3,\infty}$-norm. After  elementary
calculations, we obtain
\begin{eqnarray}
\|u(t)-w(t)\|_{3,\infty}&\leq &
\|[S(t)-S_\ell (t) S(t)]u_0\|_{3,\infty}\nonumber\\
&&+\left\|\int_0^t \left[S(t-\tau)-S_\ell(t-\tau)S(t-\tau)\right]\P 
F(\tau)\;d\tau\right\|_{3,\infty}
\nonumber\\
&&+ \left\| \int_0^t \nabla\P S(t-\tau) (u\otimes u-w\otimes
w)(\tau)\;d\tau\right\|_{3,\infty}\label{th-h-0}\\
&& + \left\| \int_0^t \nabla\P S(t-\tau) (S_\ell(t-\tau)-I)(w\otimes
w)(\tau)\;d\tau\right\|_{3,\infty}.\nonumber
\end{eqnarray}

It follows from Lemma \ref{SaS} that
\begin{eqnarray*}
\|[S_\ell(t) S(t) -S(t)]u_0\|_{3,\infty} &\leq& \|p_\ell(t)* p(t/2)
-p(t/2)\|_1 \|S(t/2) u_0\|_{3,\infty} \\
&\leq& Ct^{-(1/2-1/\ell)} \|u_0\|_{3,\infty} \to 0 \quad \mbox{as}\quad
t\to\infty.
\end{eqnarray*}

The third and the fourth term on the right-hand-side of \rf{th-h-0}
tend to 0 as $t\to\infty$ in view of Lemma \ref{Sl-f} applied either to the
function $f(\cdot,t) =V(\cdot,t)$ or to $f(\cdot,t)= (w\otimes w)(\cdot,t)$.

We deal with the second term in \rf{th-h-0} exactly in the same way as in
the proof of Theorem \ref{asymp}. Repeating the calculations from
\rf{as-2}, \rf{as-3}, and \rf{as-4} we obtain
\begin{eqnarray*}
&&\left\| \int_0^t \nabla\P S(t-\tau) (u\otimes u-w\otimes
w)(\tau)\;d\tau\right\|_{3,\infty}\\
&& \leq
4\varepsilon C \int_0^\delta (1-s)^{-1} \|u(ts)-w(ts)\|_{3,\infty} \;ds+
4\varepsilon \eta \sup_{\delta t \leq \tau \leq t}
\|u(\tau)-w(\tau)\|_{3,\infty}.
\end{eqnarray*}

Now, we define $A=\limsup_{t\to\infty} \|u(t)-w(t)\|_{3,\infty}$. To show
that $A=0$, it suffices to pass to the limit as $t\to\infty$ in inequality
\rf{th-h-0} and to repeat the reasoning given at the end of the proof of
Theorem \ref{asymp}.

The limit in \rf{Lp-h} should be shown for the $L^p$-norm replaced by the 
Mar\-cin\-kie\-wicz norm, first. Here, one should proceed as in the proof of Theorem \ref{asymp-m}. Next, the proof of \rf{Lp-h} is completed by the imbedding \rf{imbedding}.
\cbdu

\section{Asymptotic stability in abstract Banach spaces}
\setcounter{equation}{0}
\setcounter{remark}{0}

The asymptotic stability analysis described in previous sections can be
generalized to the case of more general Banach spaces. Below, we formulate
such a kind of results. We skip several details of proofs because they can
be found either in \cite{K} or in \cite{LR,Ma}.

The idea of constructing solutions and to study their large time behavior
is the following. We impose the conditions on the Banach space $E$ (cf.
Definitions \ref{functional} and \ref{adequate}) which guarantee that
our Cauchy problems have local-in-time solutions in the space
$C_w ([0,T), E)$ for some $T>0$.
Next, we show that a scaling property of $\|\cdot\|_E$
allows us to obtain, moreover,   global-in-time solutions  for
suitably small initial data.  To get such
results, we introduce a new Banach space of distributions which,
roughly speaking, is a homogeneous Besov type space modeled on $E$.
This approach allows us to get solutions for initial data less
regular than those from $E$. In this abstract setting, we also study
large-time behavior of constructed solutions.

\subsection{Definitions of spaces}

Our first two definitions are minor modifications of
\cite[Def.~7 \&~8, Sec.~8]{M96} and \cite[Def.~4.1]{LR}.

\begin{definition}\label{functional}
The Banach space $(E, \|\cdot\|_E)$ is said to be
functional and translation invariant
if the following three conditions are
satisfied:
\begin{list}{}{\itemsep1mm
\labelwidth6mm \labelsep2mm \parsep1ex\topsep2mm}

\item[i.]
$\S\subset E\subset\S'$
and the both inclusions are continuous.

\item[ii.]
either these two imbeddings have a dense range or $E$ is the dual 
space $F^*$ of a functional Banach space $F$ for which these two 
imbeddings have a dense range.

\item[iii.] The norm $\|\cdot\|_E$ on $E$ is translation invariant,
i.e.
$$
\mbox{for all}\;\;\;\; f\in E \;\;\;\;\mbox{and}\;\;\;\; y\in
\bbfR^n, \;\;\;\;
\|\tau_yf\|_E=\|f\|_E.
$$
\end{list}
\end{definition}

\begin{definition}\label{adequate}
We call the  space $(E, \|\cdot\|_E)$
adequate to the problem \rf{eq}--\rf{ini} if
\begin{list}{}{\itemsep1mm
\labelwidth6mm \labelsep2mm \parsep1ex\topsep2mm}
\item[i.]
it is
a functional translation invariant Banach space;
\item[ii.]
for all $f,g\in E$, the product
$ f\otimes g$ is well-defined as the tempered distribution, moreover,
there exist $T_0>0$  and a positive function $\omega\in L^1(0,T_0)$
such that
\begin{equation}
\|\P\nabla S(\tau)\cdot (f\otimes g)\|_E\leq
\omega(\tau) \|f\|_E\|g\|_E\label{iv:1}
\end{equation}
for  every  $f,g\in E$ and $\tau\in(0,T_0)$.
\end{list}
\end{definition}

Note that inequality \rf{iv:1} for the  space $E=L^p$
appeared already in the proof of Proposition \ref{eq:reg} (cf. \rf{Sfgp}).
Since we are interested in
an incompressible flow, we can
say that
   the Banach space $\P L^p=\{f\in L^p:
\nabla\cdot f=0\}$
is adequate to the Navier-Stokes system \rf{eq}--\rf{ini}
for every $p\in (3,\infty]$.

We refer the reader to the paper \cite{K} for  other examples of
Banach spaces adequate to \rf{eq}--\rf{ini}.
Moreover, the well-suited spaces introduced in \cite{C, CM} are
functional translation invariant Banach spaces in the sense of our
Definition \ref{functional} having some  additional properties. In
particular, they
satisfy a slightly stronger condition than  \rf{iv:1}, so they are
also adequate spaces in the
sense of Definition \ref{adequate}  (see \cite[Lem. 2.1]{CM}).
Several examples of the
well-suited (or adequate) spaces for the Navier-Stokes system 
\rf{eq}--\rf{ini} are also contained in the book \cite{LR}.

\begin{remark}
Here, it is worth of emphasizing that
if $E$ is a well-suited Banach space
(or, more generally, adequate for the problem \rf{eq}--\rf{ini})
then
for any initial datum $v_0\in
E$, $\div v_0=0$, there exists $T=T(\|v_0\|_E)$ and the unique
``mild" solution to the
Navier-Stokes equations in the space $\C([0,T); E)$.
Details are contained in \cite[Theorem 2.1]{CM}.
\end{remark}

In this paper, we use Banach spaces
with norms having
additional scaling properties. In order to state this fact more precisely,
     given $f:\bbfR^n\to\bbfR^n$, we define the rescaled function
\begin{equation}
f_\lambda(x)=f(\lambda x) \label{scal}
\end{equation}
for each $\lambda>0$. We extend this definition for all $f\in
\S'$ in the standard way.

\begin{definition}
Let $(E,\|\cdot\|_E)$
     be a Banach space, which can be imbedded continuously in  $ \S'$.
The norm $\|\cdot\|_{E}$ is said to have {\it the scaling
degree equal to $k$}, if
$
\|f_\lambda\|_{E}=\lambda^{k}\|f\|_{E}
$
for each $f\in E$ such that $f_\lambda\in E$ and
for all $\lambda>0$.
\end{definition}

It is evident that the usual norms of
the  spaces $L^p$, $L^{p,\infty}$,
$L^{p,q}$ (the Lorentz space),
${\cal M}^p_q$ (the homogeneous Morrey space)
     have  the scaling  degree equal to $-n/p$ (more details on these
     spaces can be found e.g. in \cite{K}).
On the other hand, the standard norm in the homogeneous Sobolev space
$
\dot H^s=\{f\in \S'\,:\, |\xi|^s\hat f(\xi)\in
L^2\}
$
has scaling degree $s-n/2$.

\begin{remark}
In our considerations below, we  systematically assume that the
norms of Banach spaces have the scaling degrees equal to some
$k\in (-1, 0)$. Since the space $L^p$ is  our
model example, to simplify the exposition, we shall assume that
$k=-3/p$ with $p>3$.
In this work,  Banach spaces endowed with  norms having  this
property  will be usually denoted by $E_p$.
\end{remark}

Let us  fix a Banach space $E\subset \S'$
and introduce  a new space of distributions
denoted by $\BEalf$ which,
loosely speaking, is a homogeneous Besov space modeled on
$E$. The definition we are going to introduce will be an important
tool in the next sections,
where  global-in-time solutions will be constructed
(for suitably small initial data) in $ \C([0,\infty); \BEalf)$.

\begin{definition}\label{def:BEa}
Let $\alpha \geq 0$. Given a
Banach space $E$
imbedded continuously in $\S'$, we define
$$
\BEalf = \{f\in \S'\,:\,
\|f\|_\BEalf\equiv\sup_{t>0}t^{\alpha/2}\|S(t)f\|_{E}<
\infty\}.
$$
\end{definition}

Let $E=L^p(\bbfR^n)$ for a moment. It follows immediately from
the  estimates of the heat semigroup
$$
\|S(t)f\|_{L^p(\bbfR^n)}\leq C(p,q) t^{-n(1/q-1/p)/2}
\|f\|_{L^q(\bbfR^n)}
$$
for each  $1\leq q\leq p\leq \infty$,
     that $L^q\subset \BEpa$  with $\alpha=n(1/q-1/p)$.
It is easy to obtain the analogous conclusions for the Marcinkiewicz,
Lorentz, or Morrey spaces applying appropriate estimates of the heat
semigroup mentioned in \cite[Section 3]{K}.
Moreover, for $E=L^p(\bbfR^n)$,
the norm $\|\cdot\|_\BEalf$ is equivalent to the standard  norm of
     the homogeneous Besov space  ${\dot B}_{p,\infty}^{-\alpha}$
introduced via a dyadic decomposition.

\begin{remark}
If $E$ has a norm with scaling degree $k$, then
$\|\cdot\|_{\BEalf}$ has  degree $k-\alpha$.
Indeed,
first we observe that
for any $f\in \S'$ and $\lambda >0$,
\begin{equation}
     S(t)f_\lambda=\left(S(\lambda^2t)f
\right)_\lambda.\label{e:lambda}
\end{equation}
Hence,
the scaling property of the norm on $E$ implies
$$
\|f_\lambda\|_\BEalf=
\sup_{t>0} t^{\alpha/2}\|S(t)f_\lambda\|_E=
\lambda^{k-\alpha} \sup_{\lambda^2t>0}
(\lambda^2t)^{\alpha/2}\|S(\lambda^2t)f\|_E=\lambda^{k-
\alpha}
\|f\|_\BEalf.
$$
\end{remark}

\subsection{The Navier-Stokes system}

Now,  assume that $\Ep$ is the Banach space adequate to the problem
\rf{eq}--\rf{ini} which norm has the order if scaling equal to
$-3/p$ with $p>3$.
For simplicity of the exposition, we suppose, moreover, the special 
form of the external forces
$F=\nabla\cdot V$.
It is proved in \cite[Th.~5.1]{K} that
there exists $\varepsilon >0$ such that for each $v_0\in\BEpn$
and $V(t)\in E_p$ satisfying
$$
\|v_0\|_\BEpn +\sup_{t>0}t^{1-3/(2p)}\|V(t)\|_{E_p}<\varepsilon
$$
the Cauchy problem \rf{eq}--\rf{ini} has a solution
$v(x,t)$ in the space
\begin{equation}
\X\equiv  \C([0,\infty):\BEpn)\cap
\{v:(0,\infty)\to E_p\,:\,
\sup_{t>0}t^{(1-3/p)/2}\|v(t)\|_{E_p}< \infty\}.
\label{spaceX}
\end{equation}
This is the unique solution satisfying the condition
$\sup_{t>0}t^{(1-n/p)/2}\|v(t)\|_{E_p}\leq 2 \varepsilon.$

In \cite{K},
global-in-time solutions  to the Cauchy problem
\rf{eq}--\rf{ini} are obtained using the
standard argument involving the integral equation \rf{duhamel} and
Lemma \ref{lem:xyB}.
The necessary estimate of the bilinear form
$B(\cdot,\cdot)$ are derived directly from inequality \rf{iv:1}
combined
with the scaling property of the norm in $E_p$. Details are  gathered 
in \cite{K}.

The main result on the large time behavior is contained in the 
following theorem.

\begin{theorem} {\rm (\cite[Th.~6.1]{K})} \label{asymp-w}
Let the above assumptions  remain valid.
Assume that $v$ and $\tilde v$ are two solutions of \rf{eq}--\rf{ini}
corresponding to the initial data $v_0,\tilde v_0\in\BEpn$ and forces
$F=\nabla\cdot V,\;\;  \tilde F =\nabla\cdot \tilde V$,
respectively.
Suppose that
$$
\lim_{t\to\infty}t^{(1-3/p)/2}\|S(t)(v_0-\tilde 
v_0)\|_\Ep+t^{(1-3/(2p)}\|V(t)+\tilde V(t)\|_{E_p}=0.
$$
Then
$$
\lim_{t\to\infty}t^{(1-3/p)/2}\|v(\cdot,t)-\tilde v(\cdot,
t)\|_\Ep=0.
$$
\end{theorem}

The following  lemma plays an important role in the proof of
Theorem \ref{asymp-w}.

\begin{lemma}
Let  $w\in L^1(0,1)$, $w\geq 0$, and $\int_0^1w(x)\;dx<1$.
Assume that $f$ and $g$ are two nonnegative, bounded functions such
that
\begin{equation}
f(t)\leq g(t)+\int_0^1 w(\tau) f(\tau t)\;d\tau.\label{L:0}
\end{equation}
Then $\lim_{t\to\infty}g(t)=0$ implies $\lim_{t\to\infty} f(t)=0$.
\label{f:g}
\end{lemma}

We refer the reader to \cite{K} for the elementary  proof of this
lemma. Now, to
show Theorem \ref{asymp-w},
we apply Lemma \ref{f:g} with
$f(t)=t^{(1-3/p)/2}\|v(t)-\tilde v(t)\|_{E_p}$
and
\begin{eqnarray*}
g(t)&=&
t^{(1-3/p)/2}\|S(t)(v_0-\tilde v_0)\|_\Ep\\
&& +
t^{(1-3/p)/2}
\left\|\int_0^t \P\nabla S(t-\tau)\cdot \big(V(\tau)-\tilde
V(\tau)\big)
\;d\tau\right\|_{E_p}.
\end{eqnarray*}
Here,  the estimates which appear in the proof of 
global-in-time-solutions to \rf{eq}--\rf{ini}
play again the crucial role in our reasoning.

\subsection{The mollified Navier-Stokes system}
Here, the Banach space $E$ is said to be  adequate to the mollified 
problem \rf{eq-m}--\rf{ini-m} if it satisfies all the
conditions from Definitions \ref{functional} and \ref{adequate}. If 
we recall the inequality
$$
\|\omega*f\|_E\leq \|w\|_1\|f\|_E=\|f\|_E
$$
valid for every $f\in E$, we immediately obtain global-in-time 
solutions to the mollified problem \rf{eq-m}--\rf{ini-m} in the
space $\X$ defined in \rf{spaceX}
under suitable smallness assumptions on initial conditions and 
external forces. Here, it suffices only to repeat the reasoning from 
the previous subsection.

In the analysis of  the large time asymptotics, however, we should 
impose an additional assumption on the adequate Banach
space. To prove a counterpart of Theorem~\ref{asymp-m}, we should 
guarantee that standard approximations of the Dirac delta
converge in $E_p$ (cf. \rf{a-m-4}). It is well-known that this fact 
is valid if test functions are dense in $E_p$.

\begin{theorem}
Assume that the imbedding $\S\subset E$ is dense. Let $u_0\in \BEpn$ 
and $V(t)\in E_p$ satisfy
$$
\|u_0\|_{\BEpn}+\sup_{t>0}t^{1-3/(2p)}\|V(t)\|_{E_p}<\varepsilon,
$$
where $\varepsilon$ is sufficiently small constant.
Suppose that $u_0$ is homogeneous of degree $-1$ and $F=\nabla\cdot 
V$ satisfies \rf{F-scal}.
Denote by $u(x,t)=t^{-1/2}U(x/\sqrt t)$ and $v(x,t)$  respectively 
the unique solutions to
\rf{eq}--\rf{ini} and to \rf{eq-m}--\rf{ini-m}, both
corresponding to the same initial datum $u_0$ and external force $F=\div V$.
Then
$$
\lim_{t\to\infty}t^{(1-3/p)/2}\|u(\cdot,t)-w(\cdot,
t)\|_\Ep=0.
$$\cbdu
\end{theorem}
\subsection{The Navier-Stokes system with hyperdissipation}

Recall first  that solutions to the regularized Navier-Stokes
system \rf{eq-h}--\rf{ini-h} satisfy the integral equation
\rf{duhamel-h}.
If $E$ is a functional translation
invariant
Banach space (cf. Definition \ref{functional}), we have
\begin{equation}
\|S_\ell(t) S(t) f\|_E \leq C_\ell\|S(t)f\|_E\label{SlS:E}
\end{equation}
for all $f\in E$, $t>0$, and a constant $C_\ell\geq 1$ defined in \rf{Cl}.
Hence, every Banach space adequate to the Navier-Stokes system
\rf{eq}--\rf{ini} is also adequate to the system with hyperdissipation
\rf{eq-h}--\rf{ini-h}. In other words, if inequality \rf{iv:1} holds
true for the heat semigroup $S(t)$ and a functional Banach space
$E$, it is also true for $S(\tau)$ replaced by $S_\ell(\tau)S(\tau)$
and $\omega(\tau)$ replaced by $C\omega(\tau)$.
This implies that
all estimates needed in the analysis of the Navier-Stokes system 
\rf{eq}--\rf{ini}
remain true, if we replace the heat semigroup $S(t)$ by $S_\ell(t)S(t)$. One
should remember, however, that constants in all inequalities may
increase in such a new setting.

Now, we would like to compare solutions to the models
\rf{eq}--\rf{ini} and \rf{eq-h}--\rf{ini-h} as $t\to\infty$.

\begin{theorem}
Assume that $n\ell/2<p$. Let $u_0\in \BEpn$ and $V(t)\in E_p$ satisfy
$$
\|u_0\|_{\BEpn}+\sup_{t>0}t^{1-3/(2p)}\|V(t)\|_{E_p}<\varepsilon,
$$
where $\varepsilon$ is sufficiently small constant.
Denote by $u(x,t)$ and $w(x,t)$ respectively the unique solutions to
\rf{eq}--\rf{ini} and to \rf{eq-h}--\rf{ini-h}, both
corresponding to the same initial datum $u_0$ and the external force 
$F=\div V$.
Then
$$
\lim_{t\to\infty}t^{(1-3/p)/2}\|u(\cdot,t)-w(\cdot,
t)\|_\Ep=0.
$$
\end{theorem}

Lemma \ref{SaS}, saying that  the
semigroup generated by the operator $\Delta -(-\Delta)^{\ell/2}$ can
be well-approximated in $L^1$ by the heat semigroup $S(t)$,
is again an important tool in the proof of this theorem.
Using this fact we are able to derive an integral inequality
of the form \rf{L:0} for the function
$f(t)=t^{(1-3/p)/2}\|u(\cdot,t)-w(\cdot, t)\|_\Ep$. It is important
in computations that the function $u$ and $w$ satisfy the integral
equations \rf{duhamel} and \rf{duhamel-h}, respectively. Finally,
Lemma \ref{f:g} completes the proof.


\medskip
\noindent
{\bf Acknowledgements.}~The preparation of this paper was partially
supported by the KBN
grant 2/P03A/002/24, the POLONIUM project \'EGIDE--KBN
2002--2003, and
the EU contract HYKE No. HPRN-CT-2002-00282.

\end{document}